\documentstyle[12pt]{article}
\newcommand{\qed}{\hfill\rule{4pt}{8pt}\par\vspace{\baselineskip}}

\newtheorem{lemma}{Lemma}[section]
\newtheorem{theorem}[lemma]{Theorem}
\newtheorem{proposition}[lemma]{Proposition}
\newtheorem{co}[lemma]{Corollary}
\newtheorem{definition}[lemma]{Definition}
\newtheorem{remark}[lemma]{Remark}
\newtheorem{remarks}[lemma]{Remarks}

\def\Box{\mbox{$\sqcap\!\!\!\!\sqcup$}}

\def\ot{\otimes}
\def\ra{\rightarrow}
\def\al{\alpha}
\def\eps{\varepsilon}

\def\bea{\begin{eqnarray*}}
\def\eea{\end{eqnarray*}}

\begin{document}
\title{Symmetric Coalgebras}

\author{F. Casta\~{n}o Iglesias$^1$, S.
D\u{a}sc\u{a}lescu$^2$\thanks{ On leave from University of
Bucharest, Dept. Mathematics. } and
C. N\u{a}st\u{a}sescu$^3$\\[2mm]
$^1$ Departamento de Matem\'atica Aplicada, Universidad de
Almeria,\\ E 04120 -Almeria, Spain,
e-mail: fci@ual.es\\
$^2$ Kuwait University, Faculty of Science, Dept. Mathematics,
PO BOX 5969,\\ Safat 13060, Kuwait, e-mail: sdascal@mcs.sci.kuniv.edu.kw\\
$^3$  University of Bucharest, Facultatea de Matematica, Str.
Academiei 14,\\ Bucharest 1, RO-70109, Romania, e-mail:
cnastase@al.math.unibuc.ro }

\date{}
\maketitle

\begin{abstract}
We construct a structure of a ring with local units on a
co-Frobenius coalgebra. We study a special class of co-Frobenius
coalgebras whose objects we call symmetric coalgebras. We prove
that any semiperfect coalgebra can be embedded in a symmetric
coalgebra. A dual version of Brauer's equivalence theorem is
presented, allowing a characterization of symmetric coalgebras by
comparing certain functors. We define an automorphism of the ring
with local units constructed from a co-Frobenius coalgebra, which
we call the Nakayama automorphism. This is used to give a new
characterization to symmetric coalgebras and to describe Hopf
algebras that are symmetric as coalgebras. As a corollary we
obtain as a consequence the known characterization of Hopf
algebras that are symmetric as
algebras.
\end{abstract}

\section{Introduction and Preliminaries}

Frobenius algebras appeared in group representation theory around
100 years ago. Afterwards they were recognized in many fields of
mathematics: commutative algebra, topology, quantum field theory,
von Neumann algebras, Hopf algebras, quantum Yang-Baxter equation,
etc; see \cite{cr}, \cite{lam} for classical aspects, and
\cite{cmz}, \cite{fms}, \cite{kad}, \cite{par} for more recent
developments. For instance in \cite{kad} ideas about Frobenius
algebras, Hopf subalgebras, solutions of the Yang-Baxter equation,
the Jones polynomial and 2-dimensional topological quantum field
theories are connected. An important class of algebras, which lies
between the class of Frobenius algebras and the class of
semisimple algebras, consists of symmetric algebras, as showed by
Eilenberg and Nakayama in 1955. Symmetric algebras play a special
role in representation theory by the fact that for such a $k$
algebra $R$, the $k$-dual functor is naturally equivalent to the
$R$-dual functor; we refer to the monographs \cite{cr} and
\cite{lam} and the references indicated there. For Hopf algebras,
the characterization of finite dimensional Hopf algebras that are
symmetric as algebras was given by Oberst and Schneider in
\cite{os}. The dual concept of co-Frobenius coalgebra, was
initiated by Lin in \cite{lin} for not necessarily finite
dimensional coalgebras. Several homological characterizations and
properties have been evidenced for such coalgebras. A fundamental
result, that emphasized the role of this class of coalgebras, is
that a Hopf algebra is co-Frobenius as a coalgebra if and only if
it has non-zero integrals. This was based on the previous work
\cite{ls} of Larson and Sweedler on integrals for Hopf algebras.
This fact led to a study of Hopf algebras with non-zero integrals
from a coalgebraic point of view, which has proved to be very
efficient and produced natural easy proofs of facts like the
uniqueness of the integrals. The main aim of this paper is to
define and study a special class of co-Frobenius coalgebras which
we call symmetric coalgebras. In the finite dimensional case, they
are precisely duals of symmetric algebras. In the infinite
dimensional case some completely new aspects show up. We also
study Hopf algebras that are symmetric as
coalgebras.\\

In Section 1 we review some facts about semiperfect and
co-Frobenius coalgebras, and we describe the connection to
bilinear forms. A characterization for left and right semiperfect
coalgebras that are left and right co-Frobenius is given. In
Section 2, we define two structures of a ring (without unit) on a
co-Frobenius coalgebra, by transfering the structure of a ring
with local units from the rational part of the dual algebra. We
prove that in fact the two structures that we obtain are the same.
The idea to transfer the multiplication in this way already
appeared  for compact quantum groups in \cite{pw}, where the
convolution product was defined via the inverse of the Fourier
transform. The fact was extended in \cite{an} to Hopf algebras
with non-zero integrals. A dual construction was studied in
\cite{ab}, where a coalgebra structure was constructed on a
Frobenius algebra by transporting back the coalgebra structure of
the dual coalgebra. In Section 3 we define the concept of
symmetric coalgebra. We give equivalent conditions that define
this concept, including one which uses some sort of a trace map
with respect to the ring structure of a co-Frobenius coalgebra,
and then we present several constructions that produce symmetric
coalgebras. In particular we show that cosemisimple coalgebras are
symmetric. In Section 4 we prove that any semiperfect coalgebra
can be embedded in a symmetric coalgebra by taking a certain
trivial coextension. In Section 5 we prove a result dual to
Brauer's equivalence theorem, and we give another characterization
of symmetric coalgebras by comparing some functors between
corepresentations and representations. In Section 6 we define an
automorphism of the ring with local units constructed on a
co-Frobenius coalgebra, and we call it the Nakayama automorphism.
This is used to give a new characterization of symmetric
coalgebras, and also in Section 7 to prove that a Hopf algebra $H$
is symmetric as a coalgebra if and only if it is unimodular and
the dual of the square of the antipode is an inner automorphism
with respect to the action of an invertible element of the dual
algebra $H^*$. As a consequence we obtain the result of Oberst and
Schneider that describes Hopf algebras which are symmetric as
algebras. We note that we describe explicitly the symmetric
bilinear form making $H$ a symmetric coalgebra, which might be
interesting for mathematical physicists. We also derive some
consequences concerning the antipode of a cosemisimple Hopf
algebra. Finally we discuss some facts about Hopf subalgebras that
are
symmetric as coalgebras.\\

We work over a fixed field $k$. The category of right comodules
over a coalgebra $C$ is denoted by ${\cal M}^C$. If $M$ is a right
(left) $C$-comodule, we freely regard $M$ as a left (right)
$C^*$-module, too. In particular $C$ is a $(C^*,C^*)$-bimodule,
and the left (respectively right) action of $c^*\in C^*$ on $x\in
C$ is denoted by $c^*\cdot x$ (respectively $x\cdot c^*$). A right
$C$-comodule $M$ is called quasi-finite if for any simple right
$C$-comodule $S$, the space $Hom_{C^*}(S,M)$ is finite
dimensional. If we denote the socle of $M$ by $s(M)$, this is
equivalent to the fact that any simple comodule appears with
finite multiplicity in $s(M)$. \\
If $C$ is a coalgebra and $M$ is a left (or right) $C^*$-module,
we denote by $Rat(M)$ the rational part of $M$, which is the
largest submodule of $M$ whose module structure is induced by a
right (respectively left) $C$-comodule structure. A coalgebra $C$
is called right semiperfect (see \cite{lin}) if the category
${\cal M}^C$ has enough projectives, and this is equivalent to the
fact that the injective envelope of any simple left $C$-comodule
is finite dimensional. $C$ is right semiperfect if and only if
$Rat(_{C^*}C^*)$ is dense in $C^*$ in the finite topology. A
coalgebra $C$ is called left semiperfect if the opposite coalgebra
$C^{cop}$ is right semiperfect. If $C$ is left semiperfect and
right semiperfect, we simply say that $C$ is semiperfect. It is
known (see \cite[Corollary 3.2.17]{dnr}) that if $C$ is
semiperfect, then $Rat(_{C^*}C^*)=Rat(C^*_{C^*})$, and we denote
this by $Rat(C^*)$. Also $Rat(C^*)$ is a ring with local units,
i.e. for any finite subset $X$ of $Rat(C^*)$ there exists an
idempotent $e\in Rat(C^*)$ such that $ex=xe=x$ for any $x\in X$.
These idempotent elements (local units) can be defined in terms of
the injective envelopes of the simple left (or right)
$C$-comodules. Most of the coalgebras we work with are
semiperfect. However, since we also prove results for arbitrary
coalgebras, we prefer to mention each time what sort of a
coalgebra we work with. For basic facts and notation about
coalgebras and Hopf algebras we refer to \cite{dnr} and
\cite{mont}.

\section{Co-Frobenius coalgebras}

We first need the following general result.
\begin{lemma} \label{lemaquasifinite}
Let $C$ be an arbitrary coalgebra. Let $M$ be an injective and
quasifinite right $C$-comodule, and let $N$ be a right
$C$-comodule isomorphic to $M$. If $u:M\ra N$ is an injective
morphism of right $C$-comodules, then $u$ is an isomorphism.
\end{lemma}
{\bf Proof:} We have that $s(M)\simeq u(s(M))\subseteq s(N)$.
Since $M$ is quasifinite, $N$ must be quasifinite and $s(M)\simeq
s(N)$. Thus the finite multiplicity of any simple comodule is the
same in $s(M)$ and $s(N)$, showing that $u(s(M))=s(N)$. On the
other hand, since $M$ is injective, there exists a right
$C$-comodule $Y$ such that $N=u(M)\oplus Y$. If $Y\neq 0$, then
$Y$ contains a simple subcomodule, contradicting the fact that
$u(s(M))=s(N)$. Thus $Y=0$, and then $u$ is an isomorphism. \qed

We recall that a coalgebra $C$ is called left (right) co-Frobenius
if there exists a monomorphism of left (right) $C^*$-modules from
$C$ to $C^*$. If $C$ is left and right co-Frobenius we say that
$C$ is co-Frobenius. A left (right) co-Frobenius coalgebra is left
(right) semiperfect. In particular, for any co-Frobenius coalgebra
$C$ we have that $Rat(_{C^*}C^*)=Rat(C^*_{C^*})$, which we denote
by $Rat(C^*)$. Note that if $C$ is left co-Frobenius via $\al
:C\ra C^*$, then in fact $\al$ is a morphism from $C$ to
$Rat(_{C^*}C^*)$. It is proved in \cite[Theorem 2.1]{gmn} that for
a co-Frobenius coalgebra we have that $C$ and $Rat(C^*)$ are
isomorphic as left $C^*$-modules, and also as right $C^*$-modules.
Since $C$ is injective and quasifinite as a left (or right)
$C$-comodule, an immediate consequence of Lemma
\ref{lemaquasifinite} is the following.

\begin{co}  \label{coroiso}
Let $C$ be a co-Frobenius coalgebra and $\al :C\ra C^*$ a
monomorphism of left (right) $C^*$-modules. Then $Im(\al
)=Rat(C^*)$, thus $\al$ induces an isomorphism from $C$ to
$Rat(C^*)$.
\end{co}
We also recall that a bilinear form $B:C\times C\ra k$ is called:\\
\indent $\bullet$ {\it $C^*$-balanced} if $B(x\cdot
c^*,y)=B(x,c^*\cdot y)$ for any $x,y\in C$, $c^*\in C^*$.\\
\indent $\bullet$ {\it left} (resp. {\it right}) {\it
non-degenerate} if $B(C,y)=0$ (resp. $B(y,C)=0$) implies that
$y=0$.\\
\indent $\bullet$ {\it non-degenerate} if it is left
non-degenerate and right non-degenerate.\\
\indent $\bullet$ {\it symmetric} if $B(x,y)=B(y,x)$ for any
$x,y\in C$.\\

The following two lemmas are standard results about bilinear forms
(see for example \cite[Proposition 1]{lin} or \cite[Lemma
1]{doi}). We state them for a semiperfect coalgebra $C$.

\begin{lemma}  \label{bijcor1}
There is a bijective correspondence between the morphisms of left
$C^*$-modules $\al :C\ra Rat(C^*)$ and the bilinear forms
$B:C\times C\ra k$ that are $C^*$-balanced. The correspondence is
described by $B(x,y)=\al (y)(x)$. Moreover, $\al$ is injective if
and only if $B$ is left non-degenerate.
\end{lemma}
{\bf Proof:} We see that $B(x\cdot c^*,y)=\al (y)(x\cdot
c^*)=(c^*\al (y))(x)$ and $B(x,c^*\cdot y)=\al (c^*\cdot y)(x)$
for any $x,y\in C$, $c^*\in C^*$, showing that $\al$ is left
$C^*$-linear if and only if $B$ is $C^*$-balanced. The last part
follows from the fact that $\al (y)=0$ if and only if $B(C,y)=0$.
\qed

Similarly one obtains the right hand side version of the above
result.
\begin{lemma}  \label{bijcor2}
There is a bijective correspondence between the morphisms of right
$C^*$-modules $\beta :C\ra Rat(C^*)$ and the bilinear forms
$B:C\times C\ra k$ that are $C^*$-balanced. The correspondence is
described by $B(x,y)=\beta (x)(y)$. Moreover, $\beta$ is injective
if and only if $B$ is right non-degenerate.
\end{lemma}

Let now $C$ be a co-Frobenius coalgebra, and $\al :C\ra C^*$ be an
injective morphism of left $C^*$-modules. Let $B:C\times C\ra k$
be the left non-degenerate, $C^*$-balanced bilinear form induced
by $\al$, i.e. $B(x,y)=\al (y)(x)$ for any $x,y\in C$.

\begin{lemma}
$B$ is also right non-degenerate.
\end{lemma}
{\bf Proof:} We have to show that $B(x,C)=0$ implies $x=0$.
Indeed, we then have $\al (y)(x)=0$ for any $y\in C$. Since
$Im(\al )=Rat(C^*)$ by Corollary \ref{coroiso}, and $Rat(C^*)$ is
dense in $C^*$, we see that $c^*(x)=0$ for any $c^*\in C^*$,
showing that $x=0$. \qed

Now by Lemma \ref{bijcor2}, the map $\beta :C\ra C^*$, $\beta
(x)(y)=B(x,y)=\al (y)(x)$ is an injective morphism of right
$C^*$-modules, and again by Corollary \ref{coroiso} it induces an
isomorphism between $C$ and $Rat(C^*)$.

\begin{proposition}
Let $C$ be a semiperfect coalgebra. Then $C$ is co-Frobenius if
and only if there exists a non-degenerate bilinear form $B:C\times
C\ra k$ which is $C^*$-balanced.
\end{proposition}
{\bf Proof:} If $C$ is co-Frobenius, the existence of $B$ was
proved above. For the other way around, if such a form $B$ exists,
then $C$ is left co-Frobenius by Lemma \ref{bijcor1} and right
co-Frobenius by Lemma \ref{bijcor2}. \qed

\begin{remark}
There is another way we can construct the morphism of right
$C^*$-modules $\beta:C\ra C^*$ from the morphism of left
$C^*$-modules $\al:C\ra C^*$. If we regard $\al :C\ra Rat(C^*)$ as
an isomorphism of left $C^*$-modules, then it induces an
isomorphism of right $C^*$-modules between the duals, and hence an
isomorphism of right $C^*$-modules $\al ^*:Rat(Rat(C^*)^*)\ra
Rat(C^*)$. Now since $C$ is left and right semiperfect, by
\cite[Theorem 3.5]{gn} there exists an isomorphism of right
$C^*$-modules $\sigma _C:C\ra Rat(Rat(C^*)^*)$ defined by $\sigma
_C(c)(c^*)=c^*(c)$. By composing these maps, we obtain a
monomorphism of right $C^*$-modules $\beta:C\ra Rat(C^*)$, $\beta
=\al ^*\sigma_C$. We have that $\beta (x)(y)=((\al ^*\sigma
_C)(x))(y)=(\sigma _C(x)\al )(y)=\al (y)(x)$, so we obtain exactly
the morphism $\beta$ from above.
\end{remark}

\section{A ring structure on a co-Frobenius coalgebra}

In this section we assume that $C$ is a co-Frobenius coalgebra. As
in the previous section, we denote by $\al :C\ra C^*$ an injective
morphism of left $C^*$-modules, by $B:C\times C\ra k, B(x,y)=\al
(y)(x)$ the associated non-degenerate $C^*$-balanced bilinear
form, and by $\beta :C\ra C^*, \beta (x)(y)=B(x,y)$ the associated
injective morphism of right $C^*$-modules. We can transfer to $C$
the structure of a ring without identity of $Rat(C^*)$ through the
inverses of $\al$ and $\beta$ (in fact through the isomorphisms
induced by these), obtaining two multiplications on the space $C$,
each of them making it a ring with local units. If we denote these
multiplications by $\circ$ and $\odot$, then \bea x\circ
y&=&\al ^{-1}(\al (x)\al (y))\\
&=&\al (x)\al ^{-1}(\al (y)) \;\;\; ({\rm since}\; \al \; {\rm
is}\; {\rm left}\; C^*-{\rm linear})\\
&=&\al (x)\cdot y\\
&=&\sum \al (x)(y_2)y_1 \eea and similarly \bea x\odot y&=&\beta
^{-1}(\beta (x)\beta (y))\\
&=&\beta ^{-1}(\beta (x))\cdot \beta (y)\;\;\; ({\rm since}\;
\beta \; {\rm
is}\; {\rm right}\; C^*-{\rm linear})\\
&=&x\cdot \beta (y)\\
&=&\sum \beta (y)(x_1)x_2\\
&=&\sum \al (x_1)(y)x_2 \eea

\begin{proposition}
For any $x,y\in C$ we have $x\circ y=x\odot y$, thus the two
multiplications induced on $C$ by $\al$ and $\beta$ coincide.
\end{proposition}
{\bf Proof:} Let $c^*\in C^*$ and $x,y\in C$. Then \bea c^*(x\circ
y)&=&\sum \al (x)(y_2)c^*(y_1)\\
&=&(c^*\al (x))(y)\\
&=&\al (c^*\cdot x)(y)\\
&=&B(y,c^*\cdot x) \eea and \bea c^*(x\odot y)&=&\sum \al
(x_1)(y)c^*(x_2)\\
&=&\sum \beta (y)(x_1)c^*(x_2)\\
&=&(\beta (y)c^*)(x)\\
&=&\beta (y\cdot c^*)(x)\\
&=&B(y\cdot c^*,x) \eea Since $B$ is $C^*$ balanced we see that
$c^*(x\circ y)=c^*(x\odot y)$ for any $c^*\in C^*$, which implies
that $x\circ y=x\odot y$. \qed

\begin{theorem} \label{localunits}
Let $C$ be a co-Frobenius coalgebra. Then $(C,\circ )$ is a ring
with local units such that the multiplication is a morphism of
$(C^*,C^*)$-bimodules.
\end{theorem}
{\bf Proof:} It remains to prove that $((c^*\cdot x)\circ
y)=c^*\cdot (x\circ y)$ and $x\circ (y\cdot c^*)=(x\circ y)\cdot
c^*$ for any $x,y\in C, c^*\in C^*$. For the first relation we
have that \bea
(c^*\cdot x)\circ y&=&\sum c^*(x_2)x_1\circ y\\
&=&\sum c^*(x_2)\al (x_1)(y_2)y_1\\
&=&\sum B(y_2,x_1)c^*(x_2)y_1\\
&=&\sum B(y_2,c^*\cdot x)y_1\\
&=&\sum B(y_2\cdot c^*,x)y_1\\
&=&\sum B(y_3,x)c^*(y_2)y_1\\
&=&\sum \al (x)(y_3)c^*(y_2)y_1\\
&=&c^*\cdot (\sum \al (x)(y_2)y_1)\\
&=&c^*\cdot (x\circ y) \eea while for the second one we see that
\bea x\circ (y\cdot c^*)&=&\sum c^*(y_1)x\circ y_2\\
&=&\sum c^*(y_1)\al (x)(y_3)y_2\\
&=&\sum \al (x)(y_2)y_1\cdot c^*\\
&=&(x\circ y)\cdot c^* \eea \qed

We recall that for a ring $R$ which does not necessarily have
identity, a left $R$-module $M$ is called unital if $RM=M$. The
category of unital left $R$-modules is denoted by $R-uMod$. The
following describes the unital modules for the ring $(C,\circ )$.

\begin{theorem}
Let $C$ be a co-Frobenius coalgebra. Then the category of left
unital modules $C-uMod$ associated to the ring with local units
$(C,\circ )$ is isomorphic to the category ${\cal M}^C$ of right
comodules over the coalgebra $C$.
\end{theorem}
{\bf Proof:} Since the rings $(C,\circ )$ and $Rat(C^*)$ are
isomorphic, the categories of unital left modules $C-uMod$ and
$Rat(C^*)-uMod$ are isomorphic. But for a semiperfect coalgebra
$Rat(C^*)-uMod$ is isomorphic to ${\cal M}^C$ by \cite[Proposition
2.7]{bdgn}. \qed The category isomorphism from the theorem induces
immediately the following.

\begin{co}
Let $M$ be a right $C$-comodule. Then the lattices of subobjects
of the right $C$-comodule $M$ and of the unital left $C$-module
$M$ are isomorphic.
\end{co}

\begin{co}  \label{corideale}
Let $I$ be a subspace of the co-Frobenius coalgebra $C$. Then $I$
is a left (right) ideal of the ring $(C,\circ )$ if and only if
$I$ is a right (left) coideal of the coalgebra $C$.
\end{co}

\section{Symmetric coalgebras}

A finite dimensional algebra $A$ over the field $k$ is called
symmetric if there exists a symmetric non-degenerate bilinear form
$<\; ,\;>:A\times A\ra k$ such that $<ab,c>=<a,bc>$ for any
$a,b,c\in A$. We have the following characterization of symmetric
algebras.

\begin{theorem} (\cite[Theorem 16.54]{lam})
\label{symmetricalgebra} Let $A$ be a finite dimensional algebra.
The following assertions
are equivalent.\\
(1) $A$ is symmetric.\\
(2) $A$ and $A^*$ are isomorphic as $(A,A)$-bimodules.\\
(3) There exists a $k$-linear map $f:A\ra k$ such that
$f(xy)=f(yx)$ for any $x,y\in A$, and $Ker(f)$ does not contain a
non-zero left ideal.
\end{theorem}

The aim of this section is to define the dual property for
coalgebras. Note that we do not restrict to finite dimensional
coalgebras.

\begin{definition}
A coalgebra $C$ is called symmetric if there exists an injective
morphism $\al :C\ra C^*$ of $(C^*,C^*)$-bimodules.
\end{definition}
Obviously, a symmetric coalgebra is co-Frobenius. Also, for a
finite dimensional coalgebra $C$ we have that $C$ is symmetric if
and only if the dual algebra $C^*$ is symmetric. Symmetric
coalgebras can be defined in an alternative way, as the next
result shows.

\begin{theorem}
Let $C$ be a coalgebra. The following assertions are equivalent.\\
(1) $C$ is a symmetric coalgebra.\\
(2) There exists a bilinear form $B:C\times C\ra k$ which is
symmetric, non-degenerate and $C^*$-balanced.\\
(3) $C$ is co-Frobenius, and there exists a linear map $f:C\ra k$
such that denoting by $(C,\circ )$ the ring with local units
defined
as in Theorem \ref{localunits} we have\\
\indent i) $f(x\circ y)=f(y\circ x)$ for any $x,y\in C$.\\
\indent ii) $f(c^*\cdot x)=f(x\cdot c^*)$ for any $x\in C,c^*\in
C^*$.\\
\indent iii) $Ker(f)$ does not contain a non-zero left (or right)
co-ideal of $C$.
\end{theorem}
{\bf Proof:} $(1)\Rightarrow (2)$ Let $\al :C\ra C^*$ be an
injective morphism of $(C^*,C^*)$-bimodules. Define $B:C\times
C\ra k$ by $B(x,y)=\al (y)(x)$ for any $x,y\in C$. Then by Lemma
\ref{bijcor1}, $B$ is bilinear, $C^*$-balanced and left
non-degenerate. Let us consider the multiplication $\circ$ defined
on $C$ by $x\circ y=\al ^{-1}(\al (x)\al (y))=\al (x)\cdot
y=x\cdot \al (y)$. Then
$$B(x\circ z,y)=B(x\cdot \al (z),y)=B(x,\al (z)\cdot y)=B(x,z\circ
y)$$ Note that for any $x,y\in C, c^*\in C^*$ we have \bea
B(c^*\cdot x,y)&=&\al (y)(c^*\cdot x)\\&=&(\al (y)c^*)(x)\\&=&\al
(y\cdot c^*)(x)\\&=&B(x,y\cdot c^*)\eea Hence \bea B(z\circ
x,y)&=&B(\al (z)\cdot x, y)\\
&=&B(x,y\cdot \al (z))\\
&=&B(x,y\circ z)\eea Let $x,y\in C$. Since $C$ has local units,
there exists $e\in C$ with $e\circ x=x\circ e=x$ and $e\circ
y=y\circ e=y$. Then \bea B(x,y)&=&B(x,y\circ e)\\
&=&B(x\circ y,e)\\
&=&B(y,e\circ x)\\
&=&B(y,x) \eea so $B$ is symmetric. This shows that $B$ is also
right non-degenerate.\\
$(2)\Rightarrow (1)$ Define $\al :C\ra C^*$ by $\al
(y)(x)=B(x,y)$. Then $\al$ is an injective morphism of
$C^*,C^*$-bimodules by Lemmas \ref{bijcor1} and \ref{bijcor2}. \\
$(2)\Rightarrow (3)$ Clearly $C$ is co-Frobenius. Let $\al :C\ra
C^*$ be an injective morphism of $(C^*,C^*)$-bimodules and let
$\circ$ be the multiplication induced on $C$. Define $f:C\ra k$ as
follows. Let $x\in C$. Then there exists an idempotent $e\in C$
such that $x\circ e=e\circ x=x$. We set
$f(x)=B(x,e)$. \\
We first show that $f$ is well defined. Indeed, if $e'\in C$ is
another idempotent with $e'\circ x=x\circ e'=x$, then \bea
B(x,e')&=&B(e\circ x,e')\\
&=&B(x,e'\circ e)\\
&=&B(x\circ e',e)\\
&=&B(x,e) \eea We have used the relation already proved in
$(1)\Rightarrow (2)$. \\
Let $x,y\in C$, and pick an idempotent $e\in C$ such that $x\circ
e=e\circ x=x$ and $y\circ e=e\circ y=y$. Then we have \bea
f(x\circ y)&=&B(x\circ y,e)\\
&=&B(y,e\circ x)\\
&=&B(y,x)\\
&=&B(y,x\circ e)\\
&=&B(y\circ x,e)\\
&=&f(y\circ x) \eea which proves $(i)$.\\
Let now $x\in C$ and $c^*\in C^*$. We show that $f(c^*\cdot
x)=f(x\cdot c^*)$. Pick an idempotent $e\in C$ such that $x\circ
e=e\circ x=x$, $(c^*\cdot x)\circ e=e\circ (c^*\cdot x)=c^*\cdot
x$ and $(x\cdot c^*)\circ e=e\circ (x\cdot c^*)=x\cdot c^*$. Then
we have \bea f(x\cdot c^*)&=&B(x\cdot c^*,e)\\
&=&B((x\circ e)\cdot c^*,e)\\&=&B(x\circ (e\cdot
c^*),e)\\&=&B(e\cdot c^*,e\circ x)\\&=&B(e\cdot
c^*,x)\\&=&B(e,c^*\cdot x)\\&=&B(c^*\cdot x,e)\\&=&f(c^*\cdot
x)\eea Finally, to show $(iii)$, assume that $I\subseteq Ker(f)$
for a right coideal $I$ of $C$. Then by Corollary \ref{corideale}
we have that $I$ is a left ideal in the ring $(C,\circ )$. Let $
x\in I$. For any $c\in C$ pick an idempotent $e_c\in C$ such that
$x\circ e_c=e_c\circ x=x$ and $c\circ e_c=e_c\circ c=c$. Then \bea
B(c,x)&=&B(c,x\circ e_c)\\
&=&B(c\circ x,e_c)\\
&=&f(c\circ x)\\
&=&0\eea which implies that $x=0$ by the non-degeneracy of $B$.
Thus $I=0$.\\
$(3)\Rightarrow (2)$ Define $B:C\times C\ra k$ by $B(x,y)=f(x\circ
y)$. Then clearly $B$ is bilinear and symmetric. If $B(C,x)=0$,
then $f(C\circ x)=0$, so the left coideal $C\circ x$ must be zero.
Since $C$ has local units this implies that $x=0$. Thus $B$ is
non-degenerate. We finally show that $B$ is $C^*$-balanced.
Indeed, for $x,y\in C$ and $c^*\in C^*$ we have that \bea B(x\cdot
c^*,y)&=&f((x\cdot c^*)\circ y)\\&=&f(y\circ (x\cdot
c^*))\\&=&f((y\circ x)\cdot c^*)\\&=&f(c^*\cdot (y\circ
x))\\&=&f((c^*\cdot y)\circ x)\\&=&f(x\circ (c^*\cdot
y))\\&=&B(x,c^*\cdot y)\eea and this ends the proof. \qed

\begin{remark}
We note that if $C$ is a symmetric coalgebra, then $(x\cdot
c^*)\circ y=x\circ (c^*\cdot y)$ for any $x,y\in C$, $c^*\in C^*$.
Indeed we have that $(\al (x)c^*)\al (y)=\al (x) (c^*\al (y))$.
Since $\al$ is a morphism of left and right $C^*$-modules, we
obtain $\al (x\cdot c^*)\al (y)=\al (x)\al (c^*\cdot y)$. The
desired relation follows now by applying $\al ^{-1}$, which is an
algebra morphism.
\end{remark}

\begin{proposition}\label{directsum}
Let $(C_i)_{i\in I}$ be a family of symmetric coalgebras. Then
$C=\oplus _{i\in I}C_i$ is a symmetric coalgebra.
\end{proposition}
{\bf Proof:} For any $i\in I$, let $B_i:C_i\times C_i\ra k$ be a
symmetric, non-degenerate and $C_i^*$-balanced bilinear form. We
define the bilinear form $B:C\times C\ra k$ such that the
restriction of $B$ to $C_i\times C_i$ is $B_i$ for any $i\in I$,
and $B(x,y)=0$ for any $x\in C_i, y\in C_j$ with $i\neq j$. Then
clearly $B$ is symmetric, non-degenerate and $C^*$-balanced. \qed

\begin{co}\label{cosemisimple}
A cosemisimple coalgebra is symmetric.
\end{co}
{\bf Proof:} A cosemisimple coalgebra is a direct sum of simple
coalgebras. But a simple subcoalgebra is necessarily finite
dimensional, and then it is symmetric since its dual, a matrix
algebra over a division ring, is a symmetric algebra (see
\cite[Example 16.59]{lam}). \qed

\begin{remark}
We note that Proposition \ref{directsum} provides in particular an
example of an infinite dimensional symmetric coalgebra which is
not cosemisimple.
\end{remark}

\begin{proposition} \label{tensor}
Let $C$ and $D$ be symmetric coalgebras. Then $C\otimes D$ is
symmetric.
\end{proposition}
{\bf Proof:} Let $\al :C\ra C^*$ be an injective morphism of
$(C^*,C^*)$-bimodules, and let $\beta :D\ra D^*$ be an injective
morphism of $(D^*,D^*)$-bimodules. Then it is straightforward to
check that $\al \ot \beta :C\ot D\ra C^*\ot D^*\subseteq (C\ot
D)^*$ is an injective morphism of $(C^*\ot D^*,C^*\ot
D^*)$-bimodules. Since $Im(\al \ot \beta )\subseteq Rat(C^*)\ot
Rat(D^*)$, and $Rat(C^*)\ot Rat(D^*)$ is dense in $(C\ot D)^*$, we
have that $\al \ot \beta$ is additionally a morphism of $((C\ot
D)^*,(C\ot D)^*)$-bimodules. \qed

\begin{co}
Let $C$ be a symmetric coalgebra and $n$ a positive integer. Then
the comatrix coalgebra $M^c(n,C)$ (see \cite[Section 3]{dnrv}) is
symmetric.
\end{co}
{\bf Proof:} It follows by Proposition \ref{tensor}, Corollary
\ref{cosemisimple} and the fact that $M^c(n,C)\simeq C\ot
M^c(n,k)$. \qed

We note that there exists co-Frobenius coalgebras that are not
symmetric. This can be seen by taking an example of a Frobenius
(finite dimensional) algebra that is not a symmetric algebra (see
for example \cite[Example 16.66]{lam}) and consider the dual
coalgebra. However, in the cocommutative case we have the
following.

\begin{proposition}
Let $C$ be a cocommutative coalgebra. Then $C$ is symmetric if and
only if $C$ is co-Frobenius.
\end{proposition}
{\bf Proof:} If $\al :C\ra C^*$ is an injective morphism of left
$C^*$-modules, then $\al$ is also a morphism of right
$C^*$-modules, since $C^*$ is commutative.\qed

In the finite dimensional case, we have some more
characterizations of symmetric coalgebras.

\begin{proposition}
Let $C$ be a finite dimensional coalgebra. The following
assertions are equivalent. \\
(i) $C$ is symmetric.\\
(ii) There exists a cocommutative element $c\in C$ (i.e. $\sum
c_1\ot c_2=\sum c_2\ot c_1$) which does not belong to any proper
left coideal of $C$.\\
(iii) The right (or left) $C^*$-module $C$ is cyclic, generated by
a cocommutative element.
\end{proposition}
{\bf Proof:} This follows by a direct dualization of Theorem
\ref{symmetricalgebra}, so we just sketch the proof. We have
noticed that $C$ is a symmetric coalgebra if and only if $C^*$ is
a symmetric algebra. If $f:C^*\ra k$ is a linear map such that
$f(uv)=f(vu)$ for any $u,v\in C^*$ and $Ker(f)$ not containing a
non-zero left ideal, then let $c\in C$ such that $f=i(c)$, where
$i:C\ra C^{**}$ is the natural linear isomorphism. It is easy to
see that $c$ is a cocommutative element. Since $Ker(f)$ does not
contain a non-zero left ideal, we have that $(Ker(f))^{\perp}=\{
x\in C|u(x)=0\; {\rm for}\; {\rm any}\; u\in Ker(f)\}$ is not
contained in a proper left coideal. But
$Ker(f)=Ker(i(c))=c^{\perp}=\{ u\in C^*|u(c)=0\}$, so then
$(Ker(f))^{\perp}=kc$, and we obtain the characterization $(ii)$.

To obtain $(iii)$, we find the cocommutative element $c$ as for
$(ii)$, and we note that if $f(C^*u)=0$ for some $u\in C^*$, we
must have $u=0$. Thus if $\sum c^*(c_1)u(c_2)=0$ for any $c^*\in
C^*$, then $u=0$. Therefore if $u\cdot c=0$, then $u=0$. If we
take a representation of $\Delta (c)$ such that the $c_1$'s are
linearly independent, this implies that if $u(c_2)=0$ for any
$c_2$ (in the certain representation of $\Delta (c)$), then $u=0$.
But this means that the $c_2$'s span $C$, i.e.
$c\cdot C^*=C$. \\
The converses $(ii)\Rightarrow (i)$ and $(iii)\Rightarrow (i)$
follow by reversing the above arguments. \qed

\section{An embedding theorem}

A result of Tachikawa says that any finite dimensional algebra is
isomorphic to a quotient of a symmetric algebra (see the book of
Lam \cite[page 443]{lam}). In this section we prove a dual result
about coalgebras. We do not restrict to finite dimensional
coalgebras. Let $C$ be an arbitrary coalgebra and let $M$ be a
$(C,C)$-bicomodule, i.e. $M$ is a left $C$-comodule with comodule
structure given by $m\mapsto \sum m_{(-1)}\ot m_{(0)}$, a right
$C$-comodule with structure given by $m\mapsto \sum m_{[0]}\ot
m_{[1]}$, and $$\sum m_{(-1)}\ot (m_{(0)})_{[0]}\ot
(m_{(0)})_{[1]}=\sum (m_{[0]})_{(-1)}\ot (m_{[0]})_{(0)}\ot
m_{[1]}$$ for any $m\in M$. We define a comultiplication and a
counit on the space $D=C\oplus M$ by \bea \Delta (c,m)&=&\sum
(c_1,0)\ot (c_2,0) + \sum (m_{(-1)},0)\ot (0,m_{(0)}) + \sum
(0,m_{[0]})\ot (m_{[1]},0)\\
\eps (c,m)&=&\eps (c)\eea for any $c\in C, m\in M$. These make $D$
into a coalgebra which we call the trivial coextension of $C$ and
$M$. Clearly $C$ is isomorphic to a subcoalgebra of $D$. The
definition of the trivial coextension is close to the definition
of the coalgebra
associated to a Morita-Takeuchi context in \cite{dnrv}.\\
Since $M$ is a $(C^*,C^*)$-bimodule, the dual space $M^*$ has an
induced structure of a $(C^*,C^*)$-bimodule, therefore we can
consider the trivial extension $C^*\oplus M^*$, which has an
algebra structure with the multiplication given by
$(c^*,m^*)(b^*,n^*)=(c^*b^*,c^*n^*+m^*b^*)$ for any $c^*,b^*\in
C^*$, $m^*,n^*\in M^*$. The next result follows by a
straightforward computation.

\begin{proposition}
The dual algebra of $D=C\oplus M$ is isomorphic to the trivial
extension $C^*\oplus M^*$.
\end{proposition}

Let us take now a semiperfect coalgebra $C$ and $M=Rat(C^*)$ with
the natural structure of a $(C,C)$-bicomodule induced by the
structures of a rational $(C^*,C^*)$-bimodule. Then we can form
the trivial coextension $C\oplus M=C\oplus Rat(C^*)$.

\begin{theorem} \label{simetrichopf}
Let $C$ be a semiperfect coalgebra. Then the trivial coextension
$C\oplus Rat(C^*)$ is a symmetric coalgebra. In particular any
semiperfect coalgebra can be embedded in a symmetric coalgebra.
\end{theorem}
{\bf Proof:} Denote $M=Rat(C^*)$, $D=C\oplus M=C\oplus Rat(C^*)$,
and identify $D^*$ with $C^*\oplus M^*$. Since $Rat(C^*)$ is dense
in $C^*$, the morphism of $C^*,C^*$-bimodules $\sigma :C\ra M^*$
defined by
$\sigma (c)(m)=m(c)$ for any $c\in C,m\in M$, is injective. \\
Now define the map $\al :D\ra D^*$ by $\al (c,m)=(m,\sigma (c))$
for any $c\in C,m\in M$. We first show that $\al$ is a morphism of
left $D^*$-modules. First note that since $\sigma$ is left
$C^*$-linear we have that
$$c^*\sigma (c)=\sigma (c^*\cdot c)=\sum \sigma (c^*(c_2)c_1)=\sum
c^*(c_2)\sigma (c_1)$$ for any $c\in C, c^*\in C^*$. Also for any
$n^*\in M^*, m\in M$ we have $\sum n^*(m_{(0)})\sigma
(m_{(-1)})=n^*m$. Indeed, for any $b^*\in M$ we have \bea \sum
(n^*(m_{(0)})\sigma (m_{(-1)}))(b^*)&=&\sum
n^*(m_{(0)})b^*(m_{(-1)})\\
&=&n^*(mb^*)\\&=&(n^*m)(b^*)\eea Now for any $(c^*,n^*)\in D^*$,
$(c,m)\in D$, we have that \bea \al ((c^*,n^*)(c,m))&=&\al (\sum
c^*(c_2)(c_1,0)+\sum n^*(m_{(0)})(m_{(-1)},0)+\sum
c^*(m_{[1]})(0,m_{[0]}))\\&=&\sum c^*(c_2)(0,\sigma (c_1))+\sum
n^*(m_{(0)})(0,\sigma (m_{(-1)}))+\sum c^*(m_{[1]})(m_{[0]},0)\\
&=&(0,c^*\sigma (c))+(0,n^*m)+(c^*m,0)\\
&=&(c^*m,n^*m+c^*\sigma (c))\\
&=&(c^*,n^*)(m,\sigma (c))\eea Similarly we show that $\al$ is
right $D^*$-linear. Since $\sigma$ is right $C^*$-linear we have
as above $\sum c^*(c_1)\sigma (c_2)=\sigma (c)c^*$ for any $C\in
C,c^*\in C^*$. We also see by a direct computation that $mn^*=\sum
n^*(m_{[0]})\sigma (m_{[1]})$. Hence \bea \al ((c,m)(c^*,n^*))&=&
\al (\sum c^*(c_1)(c_2,0)+\sum
c^*(m_{(-1)})(0,m_{(0)})+\sum n^*(m_{[0]})(m_{[1]},0))\\
&=&\sum c^*(c_1)(0,\sigma (c_2))+\sum
c^*(m_{(-1)})(m_{(0)},0)+\sum n^*(m_{[0]})(0,\sigma (m_{[1]}))\\
&=&(0,\sigma (c)c^*)+(mc^*,0)+(0,mn^*)\\
&=&(mc^*,mn^*+\sigma (c)c^*)\\&=&(m,\sigma (c))(c^*,n^*)\eea
Clearly $\al$ is injective, and we conclude that $C$ can be
embedded in the symmetric coalgebra $D$.\qed

\begin{remark}
We note that any subcoalgebra of a symmetric coalgebra is
semiperfect. This follows by the fact that a subcoalgebra of a
left (right) semiperfect coalgebra is also left (right)
semiperfect (see \cite[Corollary 3.2.11]{dnr}). Thus the condition
that $C$ is semiperfect is necessary for embedding $C$ in a
symmetric coalgebra.
\end{remark}

\section{A coalgebra version for the Brauer Equivalence Theorem}

Let $C$ be an arbitrary coalgebra. We consider the contravariant
functors
$$F:{\cal M}^C\longrightarrow {\rm mod}-C^*,\;\;
F(M)=Hom_k(_{C^*}M,k)$$
$$G:{\cal M}^C\longrightarrow {\rm mod}-C^*,\;\;
G(M)=Hom_{C^*}(_{C^*}M,_{C^*}C_{C^*})$$
$$H:{\cal M}^C\longrightarrow {\rm mod}-C^*,\;\;
H(M)=Hom_{C^*}(_{C^*}M,_{C^*}C^{*}_{C^*})$$ where $mod-C^*$
denotes the category of right $C^*$-modules. There is no danger of
confusion if we denote by $\leftarrow$ the right action of $C^*$
on any of $F(M),G(M),H(M)$. The following is a coalgebra version
for the Brauer equivalence theorem.

\begin{theorem}  \label{Brauer}
The functors $F$ and $G$ are naturally equivalent.
\end{theorem}
{\bf Proof:} For any $M\in {\cal M}^C$ we define $\al (M):G(M)\ra
F(M)$ by $\al (M)(f)=\eps f$ for any $f\in G(M)$. We see that $\al
(M)$ is a morphism of right $C^*$-modules since \bea (\al
(M)(f\leftarrow
c^*))(m)&=&(\eps (f\leftarrow c^*))(m)\\
&=&\eps (f(m)\cdot c^*)\\&=&\sum c^*(f(m)_2)\eps (f(m)_1)\\&=&\sum
c^*(m_1)\eps (f(m_0))\\
&=&(\eps f)(c^*\cdot m)\\
&=&(\al (M)(f)\leftarrow c^*)(m)\eea We also define $\beta
(M):F(M)\ra G(M)$ by $(\beta (M)(g))(m)=\sum g(m_0)m_1$ for any
$g\in F(M)$. We have that $\beta (M)(g)$ is left $C^*$-linear
since \bea (\beta (M)(g))(c^*\cdot m)&=&\sum g((c^*\cdot
m)_0)(c^*\cdot m)_1\\
&=&\sum c^*(m_2)g(m_0)m_1\\&=&c^*\cdot (\beta (M)(g))(m)\eea Now
we have \bea ((\beta (M)\al (M)))(f))(m)&=&\sum \eps (f(m_0))m_1\\
&=&\sum \eps (f(m)_1)f(m)_2\\&=&f(m)\eea and \bea ((\al (M)\beta
(M))(g))(m)&=&(\eps \beta (M)(g))(m)\\&=&\sum \eps (m_1)g(m_0)\\
&=&g(m)\eea showing that $\al (M)$ and $\beta (M)$ are inverse
each other. It is also easy to see that $\al$ defines a functorial
morphism. \qed

As a biproduct of the above proof, we describe the automorphisms
of the right $C$-comodule $C$. This will be used in the next
section. We denote by $U(A)$ the set of invertible elements of an
algebra $A$.
\begin{proposition}  \label{iso}
Let $C$ be a coalgebra. Then a map $f:C\ra C$ is an isomorphism of
left $C^*$-modules if and only if there exists $u\in U(C^*)$ such
that $f(c)=c\cdot u$ for any $c\in C$.
\end{proposition}
{\bf Proof:} Let us consider the isomorphism (of right
$C^*$-modules) $\beta (C):C^*\ra Hom_{C^*}(C,C)$ from the proof of
Theorem \ref{Brauer}. We have that $(\beta (C)(c^*))(c)=\sum
c^*(c_1)c_2=c\cdot c^*$ for any $c\in C, c^*\in C^*$. Then clearly
$\beta (C)(c^*d^*)=\beta (C)(d^*)\beta (C)(c^*)$, so $\beta (C)$
is an anti-isomorphism of algebras, and then the result follows by
taking the induced bijective correspondence between $U(C^*)$ and
$U(Hom_{C^*}(C,C))$.\qed

Now we can characterize symmetric coalgebras by using the above
functors.

\begin{theorem}
If $C$ is a symmetric coalgebra, then the functors $G$ and $H$ are
naturally equivalent. Conversely, if $C$ is a semiperfect
coalgebra and $G\simeq H$, then $C$ is a symmetric coalgebra.
\end{theorem}
{\bf Proof:} Assume that $C$ is symmetric. Then $C$ is semiperfect
and $C\simeq Rat(C^*)$ as $(C^*,C^*)$-bimodules. Then for any
$M\in {\cal M}^C$ we have
$$H(M)=Hom_{C^*}(_{C^*}M,_{C^*}C^{*}_{C^*})=Hom_{C^*}(_{C^*}M,Rat(C^*))\simeq
Hom_{C^*}(_{C^*}M,_{C^*}C_{C^*})=G(M)$$

Conversely, assume that $C$ is semiperfect and $G\simeq H$. Since
for any $M\in {\cal M}^C$ we have $H(M)=Hom_k(_{C^*}M,Rat(C^*))$,
we see that there exists an isomorphism
$$\al :Hom_{C^*}(-,_{C^*}C_{C^*})\ra Hom_{C^*}(-,Rat(C^*))$$
Define the morphism of left $C^*$-modules $u:C\ra Rat(C^*)$ by
$u=\al (C)(1_C)$, and $v:Rat(C^*)\ra C$ by $v=\al
^{-1}(Rat(C^*))(1_{Rat(C^*)})$. A standard argument using Yoneda's
Lemma shows that $v$ is the inverse of $u$.

We show that $u$ is a morphism of right $C^*$-modules. Let $c^*\in
C^*$, and consider the map $f:C\ra C$, $f(x)=x\cdot c^*$, which is
left $C^*$-linear. If we apply
$$Hom(f,1_C)\circ \al (C)=\al (C)\circ Hom(f,1_C)$$
to $1_C$, we obtain that $uf=\al (C)(f)$. We have that
$(uf)(x)=u(x\cdot c^*)$. On the other hand
$$(1_C\leftarrow c^*)(x)=1_C(x)\cdot c^*=x\cdot c^*=f(x)$$
so $1_C\leftarrow c^*=f$, and then \bea (\al (C)(f))(x)&=&\al
(C)((1_C\leftarrow c^*)(x)\\
&=&(\al (C)(1_C)\leftarrow c^*)(x)\\
&=&((\al (C)(1_C))(x))c^*\\
&=&u(x)c^*\eea Therefore $u(x\cdot c^*)=u(x)c^*$, which shows that
$u$ is a morphism of $(C^*,C^*)$-bimodules. \qed

\section{The Nakayama automorphism}

 Let $C$ be a co-Frobenius coalgebra and $\circ$ the multiplication defined on $C$ as in Section 2.
 Let $B:C\times C\rightarrow
 k$ be a non degenerate $C^*$-balanced bilinear form, and let
 $$\alpha : _{C^*}C\rightarrow_{C^*}C^* , \, \, \alpha(y)(x)=B(x,y)$$
 and
$$\beta : C_{C^*}\rightarrow C^*_{C^*} , \, \, \beta(y)(x)=B(y,x)$$
be the injective maps associated to $B$ as in Section 1.

Let $x\in C$.
 Since $\alpha(x)\in Rat(_{C^*}C^*)=Rat(C^*_{C^*} )$, there exists
 a unique $\sigma(x)\in C$ such that $ \alpha (x) =
 \beta(\sigma(x))$, i.e. $B(x,y)=B(\sigma (y),x)$ for any $y\in C$.
 Similarly, $ \beta(x) =
 \alpha(\tau(x))$ for some $\tau (x)\in C$. Then $\alpha(x) = \beta(\sigma(x))=
 \alpha((\tau\circ \sigma)(x))$, so $\tau\circ\sigma = 1_C$.
 Similarly $\sigma\circ \tau = 1_C$, hence $\sigma$ is
 bijective. Now, if $x,y\in C$, we have
 \bea
 \beta (\sigma (x\circ y))&=&\alpha(x\circ y)\\
 & =&\alpha(x)\alpha(y)\\
 &=& \beta(\sigma(x)) \beta(\sigma(y))\\
 &=&\beta (\sigma (x)\circ \sigma (y))
 \eea
 so $\sigma(x\circ y) = \sigma(x)\circ \sigma (y)$. Therefore
 $\sigma: C\rightarrow C$ is an automorphism of the ring
 $(C,\circ )$.  We call $\sigma$ the \emph{Nakayama automorphism} of
 $C$. Note that $\sigma$ depends on the choice of the bilinear
 form $B$. However, the Nakayama automorphism is determined up to
 the inner action of an invertible element of $C^*$, as the following result shows.
\begin{proposition}  \label{connectionnakayama}
Let $C$ be a co-Frobenius coalgebra, $B,B':C\times C\ra k$ be two
non-degenerate $C^*$-balanced bilinear forms, and $\sigma, \sigma
'$ the associated Nakayama automorphisms. Then there exists $u\in
U(C^*)$ such that $\sigma '(y)=\sigma (u^{-1}\cdot y\cdot u)$ for
any $y\in C$.
\end{proposition}
{\bf Proof:} Let  $\alpha' : _{C^*}C\rightarrow _{C^*}C^*$,
$\alpha'(y)(x) = B'(x,y)$.  By Proposition \ref{iso} there exists
$u\in U(C^*)$ such that $ (\alpha^{-1}\circ \alpha')(x)= x\cdot
u$, for any $x\in C$. Then $\alpha' (x) = \alpha(x\cdot u)$ and
\bea
 B'(x,y)&=&
\alpha'(y)(x)\\
& =& \alpha(y\cdot u)(x) \\
&=& B(x, y\cdot u) \eea Hence we have that \bea B(\sigma
'(y),x)&=&B(\sigma '(y),x\cdot u^{-1}u)\\
&=&B'(\sigma '(y),x\cdot u^{-1})\\
&=&B'(x\cdot u^{-1},y)\\
&=&B(x\cdot u^{-1},y\cdot u)\\
&=&B(x,u^{-1}\cdot y\cdot u)\\
&=&B(\sigma (u^{-1}\cdot y\cdot u),x)\eea and since $B$ is
non-degenerate we must have $\sigma '(y)=\sigma (u^{-1}\cdot
y\cdot u)$. \qed

\begin{proposition}  \label{symmetricnakayama}
Let $C$ be a co-Frobenius coalgebra with non-degenerate
$C^*$-balanced bilinear form $B:C\times C\ra k$, and let $\sigma$
be the associated Nakayama automorphism. Then $C$ is symmetric if
and only if there exists $u\in U(C^*)$ such that $\sigma
(x)=u^{-1}\cdot x\cdot u$ for any $x\in C$.
\end{proposition}
{\bf Proof:} Assume that the Nakayama automorphism is interior, so
$B(x,y)= B(y, u^{-1}xu)$ for any $x,y\in C$. We define the
bilinear map $B': C\times C \rightarrow k$, by $B'(x,y)=
B(u^{-1}\cdot x,y)$. Then \bea
B'(y,x)&=& B(u^{-1}\cdot y,x)\\
&=& B(u^{-1}\cdot x\cdot u,u^{-1}\cdot  y)\\
& = &B(u^{-1}\cdot x,y)\\
& =& B'(x,y) \eea so $B'$ is symmetric. $B'$ is clearly
$C^*$-balanced since so is $B$.

For the converse, assume that $C$ is symmetric, and let $B'$ be a
symmetric non-degenerate $C^*$-balanced bilinear form. The
Nakayama automorphism associated to $B'$ is the identity, and the
result follows by applying Proposition \ref{connectionnakayama} to
$B$ and $B'$. \qed

\section{Hopf algebras that are symmetric coalgebras}

We recall that a left (resp. right) integral on a Hopf algebra $H$
is an element $t\in H^*$ such that $h^*t=h^*(1)t$ (resp.
$th^*=h^*(1)t$) for any $h^*\in H^*$. A Hopf algebra is
co-Frobenius as a coalgebra if and only if it has non-zero left
(or right) integrals. In this case the dimension of the space of
left (resp. right) integrals is 1, and $H$ is called unimodular if
these two spaces of integrals are equal. Now we are able to
describe Hopf algebras whose underlying coalgebra structure is
symmetric.

\begin{theorem} \label{hopfsymmetric}
Let $H$ be a Hopf algebra with antipode $S$. Then $H$ is symmetric
as a coalgebra if and only if $H$ is unimodular and there exists
$u\in U(H^*)$ such that $S^2(h)=u^{-1}\cdot h\cdot u=\sum
u(h_1)u^{-1}(h_3)h_2$ for any $h\in H$.
\end{theorem}
{\bf Proof:} Assume that $H$ is symmetric. Let $\al :H\ra
Rat(H^*)$ be an isomorphism of $(H^*,H^*)$-bimodules, and let
$B:H\times H\ra k$, $B(x,y)=\al (y)(x)$, be the associated
bilinear form, which is symmetric, non-degenerate and
$H^*$-balanced.

Since $\al$ is a morphism of left $H^*$-modules, we have that
$h^*\al (1)=\al (h^*\cdot 1)=\al (h^*(1)1)=h^*(1)\al (1)$ for any
$h^*\in H^*$, so $\al (1)$ is a left integral on $H$. Similarly,
since $\al$ is a morphism of right $H^*$-modules, we have that
$\al (1)$ is also a right integral on $H$, thus $H$ is unimodular.

Let $t$ be a non-zero left and right integral on $H$. We have by
\cite[Proposition 5.5.4]{dnr} that $tS=t$. Since $H$ is
co-Frobenius, the bilinear form $D:H\times H\ra k$,
$D(x,y)=t(xS(y))$, is non-degenerate and $H^*$-balanced (see
\cite[Theorem 3]{lin} or \cite[Theorem 2]{doi}). We have that
\bea
D(x,y)&=&t(xS(y))\\
&=&tS(xS(y))\\&=&t(S^2(y)S(x))\\&=&D(S^2(y),x)\eea so the Nakayama
automorphism associated to $D$ is $S^2$. By Proposition
\ref{symmetricnakayama} we see that there exists $u\in U(H^*)$
such that $S^2(h)=u^{-1}\cdot h\cdot u$ for any $h\in H$.\\

For the converse, assume that $H$ is unimodular and there exists
an invertible $u\in H^*$ with $S^2(h)=u^{-1}\cdot h\cdot u$ for
any $h\in H$. Then let $t$ be a non-zero left and right integral,
and let $D:H\times H\ra k$ be the bilinear form associated to $t$.
Then as above the Nakayama automorphism associated to $D$ is
$S^2$, and by using the converse part of Proposition
\ref{symmetricnakayama}, we obtain that $H$ is symmetric as a
coalgebra. Note that by using the proof of Proposition
\ref{symmetricnakayama}, we see that a symmetric non-degenerate
$H^*$-balanced bilinear form is
$$B:H\times H\ra k,\;\; B(x,y)=D(u^{-1}\cdot x,y)=t((u^{-1}\cdot
x)S(y))$$ \qed

\begin{remarks}
(i) The condition that there exists an invertible element $u\in
H^*$ such that $S^2(h)=u^{-1}\cdot h\cdot u$ for any $h\in H$ is
equivalent to the fact that the map $(S^2)^*$ is an inner
automorphism of the algebra $H^*$.\\
(ii) If $H$ is a Hopf algebra that is symmetric as a coalgebra,
then $S^2$ is not necessarily a inner automorphism of $H$. Indeed,
let $H=k[SL_q(2)]={\cal O}_q(SL_2(k))$, the coordinate ring of
$SL_q(2)$, with $q$ not a root of 1 and the characteristic of $k$
different from 2. It is known that $H$ is a cosemisimple Hopf
algebra. Then $H$ is symmetric, but as it is proved in
\cite[Proposition 1.2]{bdgn}, no
power of the antipode is a inner automorphism.\\
(iii) As a consequence of Theorem \ref{simetrichopf} we obtain
many examples of infinite dimensional coalgebras that are
co-Frobenius but not symmetric. Indeed, we can take Hopf algebras
with non-zero integrals that are not unimodular. A large class of
examples like this was constructed in \cite{bdg} by taking
repeated Ore extensions of group algebras, and then factoring by
certain Hopf ideals.
\end{remarks}

As a consequence we obtain the characterization of finite
dimensional Hopf algebras that are symmetric as algebras. This was
initially proved in \cite{os}. Other proofs have been given in
\cite{lorenz}, \cite{far}; see also \cite{hum}. Note that the
result in \cite{far} is based on a formula concerning the Nakayama
automorphism of a unimodular Hopf algebra, which is superseded by
a general formula for the modular function on a Hopf algebra or
augmented Frobenius algebra in \cite{kadsto}. We recall that for a
finite dimensional Hopf algebra $H$, a left (resp. right) integral
in $H$ is an element $t\in H$ such that $ht=\eps (h)t$ (resp.
$th=\eps (h)t$) for any $h\in H$. We say that $H$ is unimodular in
the finite dimensional sense if the spaces of left and right
integrals in $H$ are equal.

\begin{co}
Let $H$ be a finite dimensional Hopf algebra. Then $H$ is
symmetric as an algebra if and only if $H$ is unimodular in the
finite dimensional sense and $S^2$ is inner.
\end{co}

\begin{co}
Let $H$ be a cosemisimple Hopf algebra. Then there exists an
invertible $u\in H^*$ such that $S^2(h)=u^{-1}\cdot h\cdot u$ for
any $h\in H$.\end{co}

The following result was proved in \cite{lar} in the particular
case of cosemisimple Hopf algebras by using character theory for
Hopf algebras.

\begin{co}
Let $H$ be a Hopf algebra which is symmetric as a coalgebra. Then
$S^2(A)\subseteq A$ for any subcoalgebra $A$ of $H$.
\end{co}
{\bf Proof:} Let $u\in U(H^*)$ such that $S^2(x)=u^{-1}\cdot
x\cdot u$ for any $x\in H$. Since $A$ is a subcoalgebra of $H$, it
is also an $(H^*,H^*)$-sub-bimodule of $H$, showing that
$S^2(A)\subseteq A$. \qed

We prove now a general result for semiperfect coalgebras. Recall
that if $X,Y$ are two subspaces of a coalgebra $C$ with
comultiplication $\Delta$, then the wedge $X\wedge Y$ is defined
by $X\wedge Y=\Delta ^{-1}(X\otimes C+C\otimes Y)$. Also $\wedge
^1X=X$, and $\wedge ^nX=(\wedge ^{n-1}X)\wedge X$ for $n\geq 2$.
\begin{proposition}
Let $C$ be a right semiperfect coalgebra and $A$ a finite
dimensional subcoalgebra of $C$. Then $A_{\infty}=\bigcup _{n\geq
1}\wedge ^nA$ is a finite dimensional subcoalgebra of $C$.
\end{proposition}
{\bf Proof:} Let us denote by ${\cal C}_A=\{ M\in {\cal M}^C|\rho
_M(M)\subseteq M\ot A\}$, where $\rho _M$ denotes the comodule
structure map of $M$. It is known that ${\cal C}_A$ is a closed
subcategory of ${\cal M}^C$ (see \cite[Theorem 2.5.5]{dnr}). Since
$A$ has finite dimension, ${\cal C}_A$ has finitely many types of
simple objects. On the other hand, ${\cal C}_{A_{\infty}}$ is the
smallest localizing subcategory of ${\cal M}^C$ that contains
${\cal C}_A$, so it has the same simple objects as ${\cal C}_A$.
In fact $M\in {\cal C}_{A_{\infty}}$ if and only if for any
$M'\subseteq M$, $M'\neq M$, the object $M/M'$ contains a simple
object of ${\cal C}_A$. Since $C$ is right semiperfect, we have
that $A_{\infty}$ is also right semiperfect. But $A_{\infty}$ has
finitely many types of simple right comodules, so it is finite
dimensional by \cite[Theorem 2.1]{bdgn}. \qed

As an immediate consequence we obtain the following result that
was proved in \cite{rad} with different methods.

\begin{co}
Let $H$ be a Hopf algebra with non-zero integrals and let $A$ be a
finite dimensional subcoalgebra. Then $A_{\infty}$ has finite
dimension.
\end{co}

\begin{co}  \label{corhinfty}
Let $H$ be a Hopf algebra with non-zero integrals. Then
$H_{\infty}=(k\cdot 1)_{\infty}$ is a Hopf subalgebra of finite
dimension.
\end{co}

\begin{co}
Let $H$ be a Hopf algebra which is symmetric as a coalgebra. Let
$K$ be a Hopf subalgebra such that either $G(K)=\{ 1\}$ (i.e. $K$
is irreducible) or $G(K)=G(H)$. Then $K$ is a symmetric coalgebra.
\end{co}
{\bf Proof:} It is known that $K$ has non-zero integrals (however
the integral on $K$ is not necessarily the restriction of the
integral on $H$ to $K$, see \cite{bdgn}), and in each of the two
situations $K$ is also unimodular (see \cite[Exercise
5.5.10]{dnr}). Let $u\in U(H^*)$ as in Theorem
\ref{hopfsymmetric}. If $i:K\ra H$ is the inclusion map, then
$v=i^*(u)\in U(K^*)$. Since the antipode $S_K$ of $K$ is the
restriction of the antipode $S$ of $H$, we have that
$S_K^2(x)=u^{-1}\cdot x\cdot u=v^{-1}\cdot x\cdot v$ for any $x\in
K$, so $K$ is symmetric. \qed

\begin{co}
Let $H$ be a Hopf algebra which is symmetric as a coalgebra. Then
the finite dimensional Hopf algebra $H_{\infty}$ from Corollary
\ref{corhinfty} is a symmetric coalgebra.
\end{co}

\begin{co}
Let $H$ be a Hopf algebra which is symmetric as a coalgebra. If
$K$ is a Hopf subalgebra of $H$ containing the coradical of $H$,
then  $K$ is also symmetric as a coalgebra.
\end{co}

\end{document}